\documentclass[12pt]{article}

\usepackage{amsmath,amsthm,amsfonts,amssymb,amscd,cite}
\usepackage{tikz}
\usepackage{color}
\usepackage{cite, url}

\newtheorem{theorem}{Theorem}[section]
\newtheorem{conjecture}[theorem]{Conjecture}

\begin{document}

\title{Intertwining local (adjacency) metric dimension with the clique number of a graph}

\author{
Ali Ghalavand $^{a,}$\thanks{Corresponding author;  email:\texttt{alighalavand@nankai.edu.cn}}
\and
Sandi Klav\v zar $^{b,c,d}$
\and 
Xueliang Li $^{a}$
}

\maketitle

\begin{center}
$^a$ Center for Combinatorics and LPMC, Nankai University, Tianjin 300071, China \\
\medskip

$^b$ Faculty of Mathematics and Physics,  University of Ljubljana, Slovenia\\
\medskip

$^c$ Institute of Mathematics, Physics and Mechanics, Ljubljana, Slovenia\\
\medskip

$^d$ Faculty of Natural Sciences and Mathematics,  University of Maribor, Slovenia\\
\end{center}

\begin{abstract}
Let \( G \) be a simple connected graph with order \( n(G) \), local metric dimension \( \dim_l(G) \), local adjacency metric dimension \( \dim_{A,l}(G) \), and clique number \( \omega(G) \), where  $G\not\cong K_{n(G)}$ and  $\omega(G)\geq3$. It is proved that \( \dim_{A,l}(G) \leq \left\lfloor \left(\frac{\omega(G) - 2}{\omega(G) - 1}\right)n(G)\right\rfloor\). Consequently, the conjecture asserting that the latter expression is an upper bound for \( \dim_l(G)\) is confirmed. It is important to note that there are infinitely many graphs that satisfy the equalities. 
\end{abstract}

\medskip\noindent
\textbf{AMS Math.\ Subj.\ Class.\ (2020)}: 05C12, 05C69

\section{Introduction}

Let \( G \) be a finite, simple, connected graph, along with two sets \( V(G) \) and \( E(G) \), representing its vertex set and edge set, respectively. We use the notations \( n(G) \) and \( \omega(G) \) to denote the order and clique number of \( G \), respectively. For two arbitrary vertices \( u \) and \( v \) in \( G \), the notation \( d_G(u, v) \) represents the length of a shortest path connecting \( u \) with \( v \). Two vertices \( u \) and \( v \) of $G$ are said to be {\em distinguished} by a vertex \( w \) or equivalently, \( w \) {\em distinguishes} \( u \) and \( v \) if the distances to \( w \) from \( u \) and \( v \) are not equal, that is, if \( d_G(u, w) \neq d_G(v, w) \).
 
In this study, we focus on examining the relationship between the local metric dimension and the local adjacency metric dimensions with the clique number of $G$. 
First, let's define the key concepts and discuss their histories. If \( S\subseteq V(G) \) and $\overline{S} = V(G)\setminus S$, then $S$ is 
\begin{itemize}
\item a {\em resolving set} if any two vertices from $\overline{S}$ are distinguished by a vertex from \( S \);
\item a {\em local resolving set} if any two adjacent vertices from $\overline{S}$ are distinguished by a vertex from \( S \);
\item an {\em adjacency resolving set} if for two vertices from $\overline{S}$ there exists a vertex from \( S \) that is adjacent to exactly one of them; 
\item a {\em local adjacency resolving set} if for any two adjacent vertices from $\overline{S}$ there exists a vertex from \( S \) that is adjacent to exactly one of them. 
\end{itemize}
The cardinality of a smallest resolving set, a smallest local resolving set, a smallest adjacency resolving set, and a smallest local adjacency resolving sets of \( G \) are respectively referred to as the {\em metric dimension} \( \dim(G) \), {\em local metric dimension} \( \dim_l(G) \), {\em adjacency metric dimension} \( \dim_A(G) \), and {\em local adjacency metric dimension} \( \dim_{A,l}(G) \). Note that \( \dim_l(G) \leq \dim(G) \) and \( \dim_{A,l}(G) \leq \dim_A(G) \), while from our perspective it is most important that 
\begin{equation}
\label{eq:the-key-one}    
\dim_l(G) \leq \dim_{A,l}(G)\,.
\end{equation}

The concept of metric dimension in graphs has a long history, originally defined independently in~\cite{13, 25}. Determining the metric dimension is known to be NP-complete for general graphs~\cite{19}, and remains NP-complete on planar graphs with maximum degree six~\cite{6}. Research in this field is extensive, largely due to the metric dimension's wide range of real-world applications, including robot navigation, image processing, privacy in social networks, and monitoring intruders in networks. The 2023 overview~\cite{survey2} of key results and applications of metric dimension includes more than 200 references.

Research on the metric dimension has led to the exploration of various related concepts. A survey \cite{survey1} that focuses on these variations cites over 200 papers. One particularly interesting variant is the local metric dimension, introduced in 2010 by Okamoto, Phinezy, and Zhang \cite{Okamoto1}. Similar to the standard metric dimension, the local metric dimension presents computational challenges \cite{9,10} and has been the subject of several studies \cite{Abrishami1, 3, 4, Ghalavand1, lal-2023, 17}, including research on the fractional local metric dimension \cite{javaid-2024}. In addition, and of utmost importance to us, Jannesari and Omoomi introduced the adjacency metric dimension~\cite{Jannesari2012}, see also~\cite{bermudo-2022, koam-2022}, while Fernau and Rodr\'{\i}guez-Vel\'{a}zquez extended the definition to the local adjacency metric dimension~\cite{9}. 

Okamoto et al.~\cite{Okamoto1} established several significant relationships between the local metric dimension and the clique number: (i) $\dim_l(G) = n(G) - 1$ if and only if $G \cong K_{n(G)}$; (ii) $\dim_l(G) = n(G) - 2$ if and only if $\omega(G) = n(G) - 1$; (iii) $\dim_l(G) = 1$ if and only if $G$ is bipartite; and (iv) $\dim_l(G) \geq \max\left \{ \lceil \log_2 \omega(G) \rceil, n(G) - 2^{n(G) - \omega(G)} \right\}$. Additionally, Abrishami et al.~\cite{Abrishami1} demonstrated that $\dim_l(G) \leq \frac{2}{5}n(G)$ when $\omega(G) = 2$ and $n(G) \geq 3$, while it has been proved in~\cite{Ghalavand1} that \( \dim_l(G) \leq \left( \frac{\omega(G) - 1}{\omega(G)} \right) n(G) \), with equality occurring only if \( G \cong K_{n(G)} \). The latter result was  conjectured in~\cite{Abrishami1}. 

The main motivation for the present paper is: 
\begin{conjecture} \label{con}{\rm \cite[Conjecture 2]{Ghalavand1}}
If $G$ is a graph with $n(G) \geq \omega(G)+1 \geq 4$, then 
\[\dim_l(G)\leq\left\lfloor \left(\frac{\omega(G) - 2}{\omega(G) - 1}\right)n(G)\right\rfloor.\]
\end{conjecture} 
The conjecture has recently been confirmed in~\cite{Ghalavand2} for all graphs $G$ with \(\omega(G) \in \{n(G)-1, n(G)-2, n(G)-3\}\). It has also established that 
\begin{itemize}
\item \(n(G)-4 \leq \dim_l(G) \leq n(G)-3\) when \(\omega(G) = n(G)-2\), and 
\item \(n(G)-8 \leq \dim_l(G) \leq n(G)-3\) when \(\omega(G) = n(G)-3\).
\end{itemize}
Additionally, the authors in~\cite{Ghalavand3} confirmed Conjecture~\ref{con} for \(\omega(G) = 3\) by proving the following result. 

\begin{theorem}\label{rth1}{\rm \cite{Ghalavand3}}
If \(G\) is a graph with \(n(G) \geq 4\) vertices and \(\omega(G) = 3\), then \(\dim_l(G) \leq \lfloor \frac{1}{2} n(G) \rfloor\).
\end{theorem}

By examining the structure of the subset \(S\) in the proof of Theorem\ref{rth1}, we can observe that \(S\) also serves as a local adjacency resolving set for \(G\), hence we can also state: 

\begin{theorem}\label{rth2}{\rm \cite{Ghalavand3}}
If \(G\) is a graph with \(n(G) \geq 4\) vertices and \(\omega(G) = 3\), then \(\dim_{A,l}(G) \leq \lfloor \frac{1}{2} n(G) \rfloor\).
\end{theorem}

Recently, Conjecture~\ref{con} was confirmed also for \(\omega(G) = 4\): 

\begin{theorem}\label{rth3}{\rm \cite{Ghalavand4}}
If \(G\) is a graph with \(n(G) \geq 5\) vertices and \(\omega(G) = 4\), then \(\dim_l(G) \leq \lfloor \frac{2}{3} n(G) \rfloor\).
\end{theorem}

In this paper, we extend the strategy used in the proof of Theorem \ref{rth3} to demonstrate the following theorem. In view of~\eqref{eq:the-key-one}, this also confirms Conjecture~\ref{con}.

\begin{theorem}\label{3th}
If \(G\) is a graph with \(n(G) \geq \omega(G) + 1 \geq 4\), then 
\[
\dim_{A,l}(G) \leq \left\lfloor \left(\frac{\omega(G) - 2}{\omega(G) - 1}\right)n(G)\right\rfloor.
\]
\end{theorem}

For any positive numbers \(t \geq 2\) and \(\omega \geq 2\), let $G_{t,\omega}$ be the graph obtained from the disjoint union of $t$ copies of $K_\omega$, by selecting one vertex in each of the $t$ copies and identify them into a single vertex. Clearly, $\omega(G_{t,\omega}) = \omega$,
$n(G_{t,\omega}) = t(\omega - 1) + 1$, and we can also observe that \(\dim_{A,l}(G_{t,\omega}) = \dim_l(G_{t,\omega}) = t(\omega -2)\). Since
$$\left\lfloor \left(\frac{\omega(G_{t,\omega}) - 2}{\omega(G_{t,\omega}) - 1}\right)n(G_{t,\omega})\right\rfloor = 
\left\lfloor \left(\frac{\omega - 2}{\omega - 1}\right)(t(\omega - 1) + 1)\right\rfloor = t(\omega -2)\,,$$
there are infinitely many graphs that attain equality in Theorem~\ref{3th} as well as in Conjecture~\ref{con}.

\section{Proof of Theorem~\ref{3th}}

Before we prove the theorem, here are some notations and definitions that we will use.

\begin{figure}[ht!]
\begin{center}
\begin{tikzpicture}
\clip(-7,0) rectangle (2,3.2);
\draw (-6,1)-- (-5,1);
\draw (-5,2)-- (-5,1);
\draw (-5,2)-- (-6,1);
\draw (-5.46,2.96)-- (-6,1);
\draw (-5.46,2.96)-- (-5,2);
\draw (-3.24,1.01)-- (-3.24,2.01);
\draw (-3.24,2.01)-- (-2.24,2.01);
\draw (-2.24,2.01)-- (-2.24,1.01);
\draw (-2.24,1.01)-- (-3.24,1.01);
\draw (-3.24,2.01)-- (-2.24,1.01);
\draw (-2.24,2.01)-- (-3.24,1.01);
\draw (-2.72,3)-- (-3.24,2.01);
\draw (-2.72,3)-- (-3.24,1.01);
\draw (-2.72,3)-- (-2.24,2.01);
\draw (-1.78,1.02)-- (-1.78,2.02);
\draw (-1.78,2.02)-- (-0.78,2.02);
\draw (-0.78,2.02)-- (-0.78,1.02);
\draw (-0.78,1.02)-- (-1.78,1.02);
\draw (-1.78,2.02)-- (-0.78,1.02);
\draw (-0.78,2.02)-- (-1.78,1.02);
\draw (-1.27,3.01)-- (-1.78,2.02);
\draw (-1.27,3.01)-- (-1.78,1.02);
\draw (-0.24,1)-- (-0.24,2);
\draw (-0.24,2)-- (0.77,2);
\draw (0.77,2)-- (0.76,1);
\draw (0.76,1)-- (-0.24,1);
\draw (-0.24,2)-- (0.76,1);
\draw (0.77,2)-- (-0.24,1);
\draw (0.29,2.99)-- (-0.24,1);
\draw (-4.67,1.01)-- (-3.67,1.01);
\draw (-3.67,2.01)-- (-3.67,1.01);
\draw (-3.67,2.01)-- (-4.67,1.01);
\draw (-4.17,2.98)-- (-4.67,1.01);
\draw (-5.77,0.93) node[anchor=north west] {$K_4^{-1}$};
\draw (-4.35,0.94) node[anchor=north west] {$K_4^{-2}$};
\draw (-2.97,0.97) node[anchor=north west] {$K_5^{-1}$};
\draw (-1.55,0.98) node[anchor=north west] {$K_5^{-2}$};
\draw (0.03,0.97) node[anchor=north west] {$K_5^{-3}$};
\begin{scriptsize}
\fill [color=black] (-6,1) circle (2.5pt);
\fill [color=black] (-5,1) circle (2.5pt);
\fill [color=black] (-5,2) circle (2.5pt);
\fill [color=black] (-5.46,2.96) circle (2.5pt);
\fill [color=black] (-3.24,1.01) circle (2.5pt);
\fill [color=black] (-3.24,2.01) circle (2.5pt);
\fill [color=black] (-2.24,2.01) circle (2.5pt);
\fill [color=black] (-2.24,1.01) circle (2.5pt);
\fill [color=black] (-2.72,3) circle (2.5pt);
\fill [color=black] (-1.78,1.02) circle (2.5pt);
\fill [color=black] (-1.78,2.02) circle (2.5pt);
\fill [color=black] (-0.78,2.02) circle (2.5pt);
\fill [color=black] (-0.78,1.02) circle (2.5pt);
\fill [color=black] (-1.27,3.01) circle (2.5pt);
\fill [color=black] (-0.24,1) circle (2.5pt);
\fill [color=black] (-0.24,2) circle (2.5pt);
\fill [color=black] (0.77,2) circle (2.5pt);
\fill [color=black] (0.76,1) circle (2.5pt);
\fill [color=black] (0.29,2.99) circle (2.5pt);
\fill [color=black] (-4.67,1.01) circle (2.5pt);
\fill [color=black] (-3.67,1.01) circle (2.5pt);
\fill [color=black] (-3.67,2.01) circle (2.5pt);
\fill [color=black] (-4.17,2.98) circle (2.5pt);
\end{scriptsize}
\end{tikzpicture}
\caption{ The graphs \( K_{4}^{-i} \) for \( i \in [2] \) and \( K_{5}^{-j} \) for \( j \in [3] \).}
\label{fig1}
\end{center}
\end{figure}
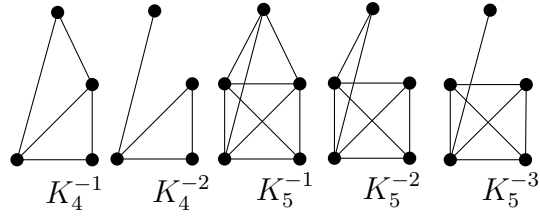

To simplify our writing, we will use the notation \( [a] \) to represent the set of natural numbers \( \{ 1, \ldots, a \} \). For clarity, we define \( [0] = \emptyset \). The notation \( d_G(u) \) refers to the degree of a vertex \( u \) in \( G \), that is, the number of vertices that are adjacent to $u$. For a subset \(W\) of the vertex set of a graph \(G\), the notation \(G[W]\) represents the subgraph of \(G\) that is induced by \(W\). 
For an arbitrary subgraph \(F\) of \(G\), the subgraph denoted \(G - F\) is defined as the graph \(G[V(G) - V(F)]\). For two vertex-disjoint subgraphs \(F_1\) and \(F_2\) of \(G\), the notation \(E_G(F_1, F_2)\) represents the set of edges in \(G\) that connect one vertex in \(F_1\) to another vertex in \(F_2\). For integers \( n \geq 4 \) and \( r \in [n-2] \), we define the graph \( K_n^{-r} \) as the graph obtained from \( K_n \) by removing $r$ edges that are incident to a common vertex. See Fig.~\ref{fig1} for some examples. 

We begin the proof by outlining a crucial approach. Let \( G \) be a graph with \( n(G)\geq \omega(G)+1\geq5 \). We are going to systematically identify maximum sets of vertex-disjoint induced subgraphs within \( G \) and its related induced subgraphs, that are isomorphic to one of the graphs \( K_{\omega(G) + 1}^{-i} \) for \( i \in [\omega(G) - 1] \) or \( K_{j} \) for \( j \in [\omega(G)] \). While this selection may not be unique, we will choose a specific selection and fix it for the purpose of this proof. This will ensure that the following notation is well-defined.

\begin{itemize}
\item Let $\mathcal{K}_1(G)$ be a maximum set of vertex disjoint induced subgraphs of $G$ isomorphic to $K_{\omega(G)+1}^{-1}$, and set $G_1 = G$.
\item For $i=2,\ldots,\omega(G)-1$, let $\mathcal{K}_i(G)$ be a maximum set of vertex disjoint induced subgraphs of $G_{i}=G_{i-1} - \sum_{X \in \mathcal{K}_{i-1}(G)}X$ isomorphic to $K_{\omega(G)+1}^{-i} $.
\item For $i=\omega(G),\ldots,2\omega(G)-1$, let $\mathcal{K}_i(G)$ be a maximum set of vertex disjoint induced subgraphs of $G_{i}=G_{i-1} - \sum_{X \in \mathcal{K}_{i-1}(G)}X$ isomorphic to $K_{2\omega(G)-i}$.
\item For $i\in[2\omega(G)-1]$, let $V_i=\cup_{X \in \mathcal{K}_i(G)}V(X)$.
\end{itemize}
Note that $\mathcal{K}_{2\omega(G)-1}(G)$ is a set of isolated vertices and that \(V_i\), \(i \in [2\omega(G)-1]\),  form a partition of \(V(G)\). It is also important to point out that some of the sets \(\mathcal{K}_i(G)\) may be empty. In the following we focus on relevant conditions concerning the elements of \(\mathcal{K}_i(G)\), for the next notation is needed. 

\begin{itemize}
  \item For $i\in[2\omega(G)-1]$, let $\mathcal{K}_i(G)=\{X^i_1,\ldots,X^i_{|\mathcal{K}_i(G)|}\}$.
  \item For indices \(i\in[\omega(G)-1]\) and \(j\in[|\mathcal{K}_i(G)|]\), define \(V(X^i_j)=\{x^i_{j_l}:l\in[|\omega(G)+1|]\}\), where the degree \(d_{X^i_j}(x^i_{j_l})\) is \(\omega(G)\) when \(l\in[\omega(G)-i]\), is \(\omega(G)-i\) when \(l=\omega(G)+1\), and equals \(\omega(G)-1\) otherwise.
    \item For indices \(i\in([2\omega(G)-1] \setminus [\omega(G)-1])\)  and \(j\in[|\mathcal{K}_i(G)|]\), define \(V(X^i_j)=\{x^i_{j_l}:l\in[2\omega(G)-i]\}\).
\end{itemize}

Let $G_i$, $i \in [2\omega(G) - 1]$, be the subgraphs of $G$ as defined above. By applying the definition of $\omega(G)$ and the maximality of $\mathcal{K}_i(G)$ for $i \in [2\omega(G) - 2]$, we can derive the following results.
\begin{enumerate}
  \item[(I)] If $i\in[|\mathcal{K}_1(G)|]$, then there is no vertex $v$ of $V(G-X^1_i)$ and a subset $X$ of $X^1_i$ such that $|X|=\omega(G)$ and $G[X\cup\{v\}]\cong\,K_{\omega(G)+1}$. 
  \item[(II)] If an integer \( \Gamma \) ranges from 2 to \( \omega(G) - 1 \), with \( i\in[|\mathcal{K}_\Gamma(G)|] \), \( j\in\{\omega(G)-\Gamma,\omega(G)-\Gamma+1\} \), and \( v\in V(G_\Gamma-X^\Gamma_i) \), then 
      \[ G[\{x_{i_l}^\Gamma:l\in([\omega(G)]-\{\omega(G)-\Gamma,\omega(G)-\Gamma+1\})\}\cup\{x_{i_j}^\Gamma,v\}]\not\cong\,K_{\omega(G)}.\]
\item [(III)] If an integer \( \Gamma \) ranges from \( \omega(G)+1 \) to \( 2\omega(G) - 1 \), and another integer \( \Upsilon \) ranges from \( \Gamma \) to \( 2\omega(G) - 1 \) with \( i \in [|\mathcal{K}_\Gamma(G)|] \), \( j \in [|\mathcal{K}_\Upsilon(G)|] \), and \( X^\Gamma_i \neq X^\Upsilon_j \), then for \( l \in [2\omega(G) - \Upsilon] \), it follows that 
\[
G[V(X^\Gamma_i) \cup \{ x^\Upsilon_{j_l} \}] \not\cong K_{\Gamma + 1}.
\]
\item[(IV)] If an integer \( \Gamma \) ranges from \( \omega(G) \) to \( 2\omega(G) - 1 \),  with \( i \in [|\mathcal{K}_{\omega(G)}(G)|] \), \( j \in [|\mathcal{K}_\Gamma(G)|] \), and \( X^{{\omega(G)}}_i \neq X^\Upsilon_j \), then for \( l_1 \in [\omega(G)] \) and \( l_2 \in [2\omega(G) - \Gamma] \), it follows that 
      \[x^{\omega(G)}_{i_{l_1}}x^{\Gamma}_{j_{l_2}}\not\in\,E(G).\]
      
  \item[(V)]  If an integer \( \Gamma \) ranges from 1 to \( \omega(G) - 1 \), where \( i\in[|\mathcal{K}_{\Gamma}(G)|]\), \( j\in[|\mathcal{K}_{\omega(G)}(G)|]\), \( l_1 \in[\omega(G)+1]\), and \( l_2 \in[\omega(G)]\), then the condition \( x^\Gamma_{i_{l_1}} \ne x^{\omega(G)}_{j_{l_2}} \) implies that there exists \( l_3 \in[\omega(G)]\) such that 
\[
G[x^\Gamma_{i_{l_1}}, x^{\omega(G)}_{j_{l_2}}, x^{\omega(G)}_{j_{l_3}}] \cong P_3.
\]
\end{enumerate}             
                            
In the following processes, we will construct a set \( S \) with \( |S| \leq \left(\frac{\omega(G)-2}{\omega(G)-1}\right)n(G) \), ensuring that \( S \) remains a  local adjacency resolving set for \( G \). We start with \( S = \emptyset \). To move forward, we need to introduce an additional notation. Let \( \mathcal{Y} \subseteq\mathcal{K}_{\omega(G)}(G)  \) and  \( F \in \mathcal{K}_{i}(G) \), where \( i\in [\omega(G)-1] \). The notation \( \tau_{i,\omega(G)}(F,\mathcal{Y}) \) represents the set of elements \( Y \) in \( \mathcal{Y} \) such that \( E_G(F,Y) \neq \emptyset \).
 
\medskip\noindent
\underline{$1^{\rm st}$ process:}
\begin{enumerate}
\item[(1.1)] Set $S = \emptyset$, $i=1$, and $\mathcal{Y} = \mathcal{K}_{\omega(G)}(G)$.
\item[(1.2)] If $i>|\mathcal{K}_{1}(G) |$, then return  $S$ and $\mathcal{Y}$,  and end the process,  otherwise go to (1.3).  
\item[(1.3)] If  $|\tau_{1,\omega(G)}(X_i^1,\mathcal{Y})|=0$, then set
\begin{align*}
S & = S\cup\,(V(X_i^1)-\{x^1_{i_{\omega(G)}},x^1_{\omega(G)+1}\}),\\
i & =i+1,
\end{align*}
and proceed to (1.2), otherwise go to (1.4).
\item[(1.4)] If  $|\tau_{1,\omega(G)}(X_i^1,\mathcal{Y})|=1$ and $\tau_{1,\omega(G)}(X_i^1,\mathcal{Y})=\{Y_1\}$, then choose two distinct elements $l_1,l_2\in[\omega(G)+1]$ and two distinct elements $y^1_1,y^1_2\in\,V(Y_1)$ such that $G[\{x^1_{i_{l_1}},y^1_1,y^1_2\}]\cong\,P_3$, and then set
\begin{align*}
S & = S\cup\,(V(X_i^1)-\{x^1_{i_{l_2}}\})\cup(V(Y_1)-\{y^1_1,y^1_2\}),\\
i & =i+1,\\
\mathcal{Y}&=\mathcal{Y}-\tau_{1,\omega(G)}(X_i^1,\mathcal{Y}),
\end{align*}
and proceed to (1.2), otherwise go to (1.5).
\item[(1.5)] If $|\tau_{1,\omega(G)}(X_i^1,\mathcal{Y})|=2$ and $\tau_{1,\omega(G)}(X_i^1,\mathcal{Y})=\{Y_1,Y_2\}$, then for an element $l\in[\omega(G)+1]$, for which $G[(V(X_i^1)-\{x^1_{i_l}\})\cup V(Y_1)\cup V(Y_2)]$ is connected, choose elements $y^1_1,y^1_2\in V(Y_1)$ and $y^2_1,y^2_2\in V(Y_2)$ such that for $z_1,z_2\in([\omega(G)+1]-\{l\})$, not necessarily distinct, we have $G[\{x^1_{i_{z_1}},y^1_1,y^1_2\}]\cong\,P_3$ and $G[\{x^1_{i_{z_2}},y^2_1,y^2_2\}]\cong\,P_3$, and then set
\begin{align*}
S & = S\cup(V(X_i^1)-\{x^1_{i_l}\})\cup(V(Y_1)-\{y^1_1,y^1_2\})\cup(V(Y_2)-\{y^2_1,y^2_2\}),\\
i & =i+1,\\
\mathcal{Y}&=\mathcal{Y}-\tau_{1,\omega(G)}(X_i^1,\mathcal{Y}),
\end{align*}
and proceed to (1.2), otherwise go to (1.6).
\item[(1.6)] If  $|\tau_{1,\omega(G)}(X_i^1,\mathcal{Y})|\geq3$ and $\tau_{1,\omega(G)}(X_i^1,\mathcal{Y})=\{Y_1,\ldots,Y_{|\tau_{1,\omega(G)}(X_i^1,\mathcal{Y})|}\}$, then for $l$ from $1$ to $|\tau_{1,\omega(G)}(X_i^1,\mathcal{Y})|$, choose two distinct elements $y^l_1,y^l_2 \in V(Y_l)$ such that for an element $z \in [\omega(G)+1]$, $G[\{x^1_{i_z},y^l_1,y^l_2\}]\cong\,P_3$, and then set
\begin{align*}
S & = S\cup\,V(X_i^1)\cup_{l=1}^{|\tau_{1,\omega(G)}(X_i^1,\mathcal{Y})|}(V(Y_l)-\{y^l_1,y^l_2\}),\\
i & =i+1,\\
\mathcal{Y}&=\mathcal{Y}-\tau_{1,\omega(G)}(X_i^1,\mathcal{Y}),
\end{align*}
and proceed to (1.2).
\end{enumerate}

For integers \( \Gamma \) ranging from 2 to \( \omega(G) - 1 \), perform the \( \Gamma^{\text{th}} \) process as follows.

\medskip\noindent
\underline{$\Gamma^{\rm th}$ process:} 
\begin{enumerate} 
\item[($\Gamma$.1)] Consider the sets $S$ and $\mathcal{Y}$ returned in the $(\Gamma-1)^{\rm th}$ process, and set $i=1$.
\item[($\Gamma$.2)] If $i>|\mathcal{K}_\Gamma(G)|$, then return  $S$ and $\mathcal{Y}$,  and end the process,  otherwise go to ($\Gamma$.3).  
\item[($\Gamma$.3)] If  $|\tau_{\Gamma,\omega(G)}(X_i^\Gamma,\mathcal{Y})|=0$, then set
\begin{align*}
S & = S\cup\,(V(X_i^\Gamma)-\{x^\Gamma_{i_{\omega(G)-\Gamma}},h^\Gamma_{i_{\omega(G)-\Gamma+1}}\}),\\
i & =i+1,
\end{align*}
and proceed to ($\Gamma$.2), otherwise go to ($\Gamma$.4).
\item[($\Gamma$.4)] If \( |\tau_{\Gamma,\omega(G)}(X_i^\Gamma,\mathcal{Y})| = 1 \) and \( \tau_{\Gamma,\omega(G)}(X_i^\Gamma,\mathcal{Y}) = \{Y_1\} \), then choose a member \( l_1 \in [\omega(G)] \) such that \(E_G( G[V(X_i^\Gamma) - \{x^\Gamma_{i_{l_1}}\}],Y_1)\neq\emptyset \). Additionally, select two distinct members \( y^1_1, y^1_2 \in V(Y_1) \) such that for a member \( l_2 \in ([\omega(G)+1] - \{l_1\}) \), it holds that
     $G[\{x^\Gamma_{i_{l_2}},y^1_1,y^1_2\}]\cong\,P_3$. Then set
\begin{align*}
S & = S\cup\,(V(X_i^\Gamma)-\{x^\Gamma_{i_{l_1}}\})\cup(V(Y_1)-\{y^1_1,y^1_2\}),\\
i & =i+1,\\
\mathcal{Y}&=\mathcal{Y}-\tau_{\Gamma,\omega(G)}(X_i^\Gamma,\mathcal{Y}),
\end{align*}
and proceed to ($\Gamma$.2), otherwise go to ($\Gamma$.5).
\item[($\Gamma$.5)] If $|\tau_{\Gamma,\omega(G)}(X_i^\Gamma,\mathcal{Y})|=2$ and $\tau_{\Gamma,\omega(G)}(X_i^\Gamma,\mathcal{Y})=\{Y_1,Y_2\}$, then for an element $l_1\in[\omega(G)]$, for which $E_G(G[V(X_i^\Gamma)-\{x^\Gamma_{i_{l_1}}\}],Y_1)\neq\emptyset$ and $E_G(G[V(X_i^\Gamma)-\{x^\Gamma_{i_{l_1}}\}],Y_2)\neq\emptyset$, choose members $y^1_1,y^1_2\in V(Y_1)$ and $y^2_1,y^2_2\in V(Y_2)$ such that for $l_2,l_3\in([\omega(G)+1]-\{l_1\})$, not necessarily distinct, we have $G[\{x^\Gamma_{i_{l_2}},y^1_1,y^1_2\}]\cong\,P_3$ and $G[\{x^\Gamma_{i_{l_3}},y^2_1,y^2_2\}]\cong\,P_3$, and then set
\begin{align*}
S & = S\cup(V(X_i^\Gamma)-\{x^\Gamma_{i_{l_1}}\})\cup(V(Y_1)-\{y^1_1,y^1_2\})\cup(V(Y_2)-\{y^2_1,y^2_2\}),\\
i & =i+1,\\
\mathcal{Y}&=\mathcal{Y}-\tau_{\Gamma,\omega(G)}(X_i^\Gamma,\mathcal{Y}),
\end{align*}
and proceed to ($\Gamma$.2), otherwise go to ($\Gamma$.6).
\item[($\Gamma$.6)] If  $|\tau_{\Gamma,\omega(G)}(X_i^\Gamma,\mathcal{Y})|\geq3$ and $\tau_{\Gamma,\omega(G)}(X_i^\Gamma,\mathcal{Y})=\{Y_1,\ldots,Y_{|\tau_{\Gamma,\omega(G)}(X_i^\Gamma,\mathcal{Y})|}\}$, then for $l$ from $1$ to $|\tau_{\Gamma,\omega(G)}(X_i^\Gamma,\mathcal{Y})|$, choose two distinct members $y^l_1,y^l_2 \in V(Y_l)$ such that for an member $z \in [\omega(G)]$, $G[\{x^\Gamma_{i_z},y^l_1,y^l_2\}] \cong\,P_3$, and then set
\begin{align*}
S & = S\cup\,V(X_i^\Gamma)\cup_{l=1}^{|\tau_{\Gamma,\omega(G)}(X_i^\Gamma,\mathcal{Y})|}(V(Y_l)-\{y^l_1,y^l_2\}),\\
i & =i+1,\\
\mathcal{Y}&=\mathcal{Y}-\tau_{\Gamma,\omega(G)}(X_i^\Gamma,\mathcal{Y}),
\end{align*}
and proceed to ($\Gamma$.2).
\end{enumerate}
Additionally, for integers \( \Gamma \) ranging from \( \omega(G) \) to \( 2\omega(G) - 3 \), execute the \( \Gamma^{\text{th}} \) process as follows.

\medskip\noindent
\underline{$\Gamma^{\rm th}$ process:} 
\begin{enumerate} 
\item[($\Gamma$.1)] Consider the set \( S \) obtained in the \((\Gamma-1)^{\rm th}\) process, and set \( i = 1 \).
\item[($\Gamma$.2)] If $i>|\mathcal{K}_{\Gamma(G)+1}|$ and $\Gamma=2\omega(G) - 3$, then set \( S = S \cup \mathcal{K}_{2\omega(G)-1}(G) \), return \( S \), and end the process. Otherwise, proceed to ($\Gamma$.3).
\item[($\Gamma$.3)] If $i>|\mathcal{K}_{\Gamma(G)+1}|$, then return  $S$,  and end the process,  otherwise go to ($\Gamma$.4).  
\item[($\Gamma$.4)] Set 
\begin{align*}
S & = S\cup(V(X^{\Gamma+1}_i)-\{x^{\Gamma+1}_{i_1}\}),\\
i & =i+1,
\end{align*}
and proceed to ($\Gamma$.2).
\end{enumerate}

Now, let \( S \) be the subset of \( V(G) \) formed from the \( 2\omega(G)-3 \) processes described above. For any positive integer \( t \), we  define \( \xi_t \) as follows:
\[
\xi_t = \frac{\omega(G)-2}{\omega(G)-1}\left(\omega(G)(t+1)+1\right).
\]
Given that \( \omega(G) > 3 \), we can derive that:
\[\xi_t\geq\left\{\begin{array}{ll}
\omega(G) - 1; & t = 0,\\
2\omega(G) - 2; & t = 1,\\
3\omega(G) - 4; & t = 2,\\
\omega(G)(t+1) - 2t + 1; & {\rm otherwise}.
                  \end{array}\right.\]
Additionally, for \( r \in ([\omega(G)-1] - \{1\}) \), we have 
\[
\frac{\omega(G)-2}{\omega(G)-1}\, r \geq r - 1.
\]
Therefore, by utilizing (I)-(V)  and applying the methods used to construct \( S \), we can confirm that \( S \) is a local adjacency resolving set for \( G \). Furthermore, 
\[
|S| \leq \frac{\omega(G) - 2}{\omega(G) - 1} n(G).
\]
By utilizing Theorem~\ref{rth2} and because \( \dim_{A,l}(G) \) is an integer, Theorem~\ref{3th} is proved. 

\section*{Acknowledgments}

The research of Ali Ghalavand and Xueliang Li was supported by the NSFC No.\ 12131013. Sandi Klav\v{z}ar was supported by the Slovenian Research Agency (ARIS) under the grants P1-0297, N1-0355, and N1-0285.

\section*{Conflicts of interest} 

The authors declare no conflict of interest.

\section*{Data availability} 

No data was used in this investigation.


\begin{thebibliography}{99}

\bibitem{Abrishami1} 
G. Abrishami, M. A. Henning, M. Tavakoli,  
Local metric dimension for graphs with small clique numbers,
Discrete Math.\ 345 (2022) Paper 112763.

\bibitem{3} 
G.A. Barrag\'{a}n-Ram\'{\i}rez, A. Estrada-Moreno, Y. Ram\'{\i}rez-Cruz, J.A. Rodr\'{\i}guez-Vel\'{a}zquez, The local metric dimension of the lexicographic product of graphs, 
Bull.\ Malays.\ Math.\ Sci.\ Soc.\ 42 (2019) 2481--2496.

\bibitem{4} 
G.A. Barrag\'{a}n-Ram\'{\i}rez, J.A. Rodr\'{\i}guez-Vel\'{a}zquez, 
The local metric dimension of strong product graphs, 
Graphs Combin.\ 32 (2016) 1263--1278.

\bibitem{bermudo-2022}
S.~Bermudo, J.M.~Rodr\'iguez, J.A.~Rodr\'iguez-Vel\'azquez, J.M.~Sigarreta, 
The adjacency dimension of graphs,
Ars Math.\ Contemp.\ 22 (2022) Paper 2.

\bibitem{6} 
J. Diaz, O. Pottonen, M. Serna, E.J. van Leeuwen, 
Complexity of metric dimension on planar graphs, 
J.\ Comput.\ Syst.\ Sci.\ 83 (2017) 132--158.

\bibitem{9} 
H. Fernau, J.A. Rodr\'{\i}guez-Vel\'{a}zquez, 
On the (adjacency) metric dimension of corona and strong product graphs and their local variants, combinatorial and computational results, 
Discrete Appl.\ Math. 236 (2018) 183--202.

\bibitem{10} 
H. Fernau, J.A. Rodr\'{\i}guez-Vel\'{a}zquez, 
Notions of metric dimension of corona products: combinatorial and computational results, 
Lecture Notes Comput.\ Sci.\ 8476 (2014) 153--166.

\bibitem{fitriani} 
D.~Fitriani, S.W.~Saputro, 
The local metric dimension of amalgamation of graphs,
Electron.\ J.\ Graph Theory Appl.\ (EJGTA) 12 (2024) 125--146.

\bibitem{Ghalavand1} 
A. Ghalavand, M. A. Henning, M. Tavakoli, 
On a conjecture about the local metric dimension of graphs,
Graphs Combin.\ 39 (2023) Paper 5.

\bibitem{Ghalavand2} 
A. Ghalavand, S. Klav\v zar, X. Li, 
Interplay between the local metric dimension and the clique number of a graph, 
\url{arXiv:2412.17074} [math.CO] (22 Dec 2024). 

\bibitem{Ghalavand3} 
A. Ghalavand, S. Klav\v zar, X. Li, 
On the local metric dimension of $K_4$-free graphs, 
\url{arXiv:2506.00414} [math.CO] (31 May 2025). 

\bibitem{Ghalavand4} 
A. Ghalavand,  X. Li, 
On the local metric dimension of $K_5$-free graphs, submitted. 

\bibitem{13} 
F. Harary, R.A. Melter, 
The metric dimension of a graph, 
Ars Combin.\ 2 (1976) 191--195.

\bibitem{Jannesari2012} M. Jannesari, B. Omoomi, 
The metric dimension of the lexicographic product of graphs, 
Discrete Math.\ 312 (2012) 3349--3356.

\bibitem{javaid-2024}
I. Javaid, H. Benish, M. Murtaza, 
The fractional local metric dimension of graphs,
Contrib.\ Discrete Math.\ 19 (2024) 163--177.

\bibitem{19} 
S. Khuller, B. Raghavachari, A. Rosenfeld, 
Landmarks in graphs, 
Discrete Appl.\ Math.\ 70 (1996) 217--229.

\bibitem{klavzar-2023} 
S. Klav\v{z}ar, D. Kuziak,
Nonlocal metric dimension of graphs,
Bull.\ Malays.\ Math.\ Sci.\ Soc.\ 46 (2023) Paper 66.

\bibitem{17} 
S. Klav\v{z}ar, M. Tavakoli, 
Local metric dimension of graphs: generalized hierarchical products and some applications, 
Appl.\ Math.\ Comput.\ 364 (2020) Paper 124676.

\bibitem{koam-2022}
A.N.A.~Koam, A.~Ahmad, M.~Azeem, A.~Khalil, M.F.~Nadeem, 
On adjacency metric dimension of some families of graph,
J.\ Funct.\ Spaces (2022) Paper 6906316.

\bibitem{survey1}
D.~Kuziak, I.G.~Yero,
Metric dimension related parameters in graphs: A survey on combinatorial, computational and applied results,  
\url{arXiv:2107.04877} [math.CO] (10 Jul 2021). 

\bibitem{lal-2023}
S. Lal, V.K. Bhat, 
On the local metric dimension of generalized wheel graph,
Asian-Eur.\ J.\ Math.\ 16 (2023) Paper 2350194.

\bibitem{Okamoto1} 
F. Okamoto, B. Phinezy, P. Zhang,  
The local metric dimension of a graph, 
Math.\ Bohem.\ 135 (2010) 239--255.

\bibitem{rodriguez-2016}
J.A.~Rodr\'iguez-Vel\'azquez, G.A.~Barrag\'an-Ram\'irez, C.~Garc\'ia G\'omez, 
On the local metric dimension of {C}orona product graphs,
Bull.\ Malays.\ Math.\ Sci.\ Soc.\ 39 (2016) S157--S173.

\bibitem{25} 
P.J. Slater, 
Leaves of trees, 
Congress.\ Numer.\ 14 (1975) 549--559. 

\bibitem{survey2}
R.~C.~Tillquist, R.~M.~Frongillo, M.~E.~Lladser. 
Getting the lay of the land in discrete space: a survey of metric dimension and its applications. 
SIAM Rev.\ 65 (2023) 919--962.

\end{thebibliography}
\end{document}